\documentclass[letterpaper,ngerman,english]{article}
\usepackage[T1]{fontenc}
\usepackage[latin9]{inputenc}
\usepackage{float}
\usepackage{mathtools}
\usepackage{amssymb}
\usepackage{graphicx}

\makeatletter

\providecommand{\tabularnewline}{\\}

\newcommand{\lyxaddress}[1]{
\par {\raggedright #1
\vspace{1.4em}
\noindent\par}
}

\makeatother

\usepackage{babel}
\begin{document}

\title{Improved Time of Arrival measurement model for non-convex optimization}

\author{Juri Sidorenko, Leo Doktorski, \\
Volker Schatz, Norbert Scherer-Negenborn, Michael Arens}
\maketitle

\lyxaddress{Fraunhofer IOSB, Ettlingen Germany \\
juri.sidorenko@iosb.fraunhofer.de, Tel.: +49 7243 992-351}


\section{Abstract}

The quadratic system provided by the Time of Arrival technique can
be solved analytically or by optimization algorithms. In practice,
a combination of both methods is used. An important problem in quadratic
optimization is the possible convergence to a local minimum, instead
of the global minimum. This article presents an approach how this
risk can be significantly reduced. The main idea of our approach is
to transform the local minimum to a saddle point, by increasing the
number of dimensions. In contrast to similar methods such as, dimension
lifting does our problem remains non-convex.

\section{Introduction}

In position estimation the Time of Arrival (ToA) \cite{TOA} technique
is standard. The area of application extends from satellite based
systems like GPS \cite{GPS1}, GLONASS \cite{Glonass}, Galileo \cite{Galileo},
mobile phone localization (GSM) \cite{GSM}, radar based systems such
as UWB \cite{UWB}, FMCW radar \cite{Radar} to acoustic systems \cite{Accustic}.
\\
The ToA technique leads to a quadratic equation. Optimization algorithms
used to solve this system depends on the initial estimate. Unfortunately
chosen initial estimates can increase the probability to convergence
to a local minimum. In some cases it is possible to transform the
quadratic to a linear system \cite{ION_JS,Sidorenko2016,Hmam2010}.
This linear system can be used to provide an initial estimate. On
the other hand, the linear system is more affected by noise, compared
to the quadratic system \cite{ION_JS,Sidorenko2016}. In practice,
a combination of both methods is used to obtain the unknown position
of an object. \cite{Abel1991,Bancroft1985,Chaffee1994}. However,
the initial estimates by a linear solution only applies if the base
station positions are known. This article presents an approach how
the risk of convergence to a local minimum during the optimization
process can be significantly reduced for the ToA technique. The approach
does not require initial estimations provided by a linear solution,
rather the insertion of an additional variable is used to transform
a local minimum to a saddle point at the same coordinates. In order
to simplify the prove of our approach, it is assumed that the position
of the base stations are known.

Our approach was inspired by dimension lifting \cite{Balas2005,Lifting,Lift}
and concave programming \cite{Liu1995}. Dimension lifting introduces
an additional dimension to transform a non-convex to a convex feasible
region. Concave programming describes a non-convex problem in terms
of d.c. functions (differences of convex functions). In our method,
the non-convex problem remains non-convex. The objective of this paper
is to show how this approach, reduces the risk of convergence to local
minimum.\\

This paper is organized as follows. The third section, introduces
the objective functions $F$ and the corresponding improved objective
functions $F_{L}$. In Section four, we use Levenberg-Marquardt algorithm
\cite{Levenberg-Marquardt} to illustrate the optimization steps for
$F$ and $F_{L}$. The last section address the results of the optimization
algorithm with randomly selected constellations. 

\section{Methodology}

\begin{table}[H]
\begin{centering}
\caption{Used notations\label{tab:Exmplation-of-the}}
\ \\
\par\end{centering}
\centering{}%
\begin{tabular}{|c|c|}
\hline 
Notations & \selectlanguage{ngerman}%
Definition\selectlanguage{english}%
\tabularnewline
\hline 
\hline 
$x,y,z$ & Estimated position of object $T$ \tabularnewline
\hline 
$x_{G},y_{G},z_{G}$ & Ground truth position of object $T$ \tabularnewline
\hline 
$a_{i},b_{i},c_{i}$ & Ground truth position of base stations $B_{i}$, $1\leq i\leq N$\tabularnewline
\hline 
$d_{i}$ & Distance measurements between base stations $B_{i}$ and object $T$\tabularnewline
\hline 
$\lambda$ & Additional variable\tabularnewline
\hline 
\end{tabular}
\end{table}

\begin{figure}[H]
\centering{}\includegraphics[scale=0.5]{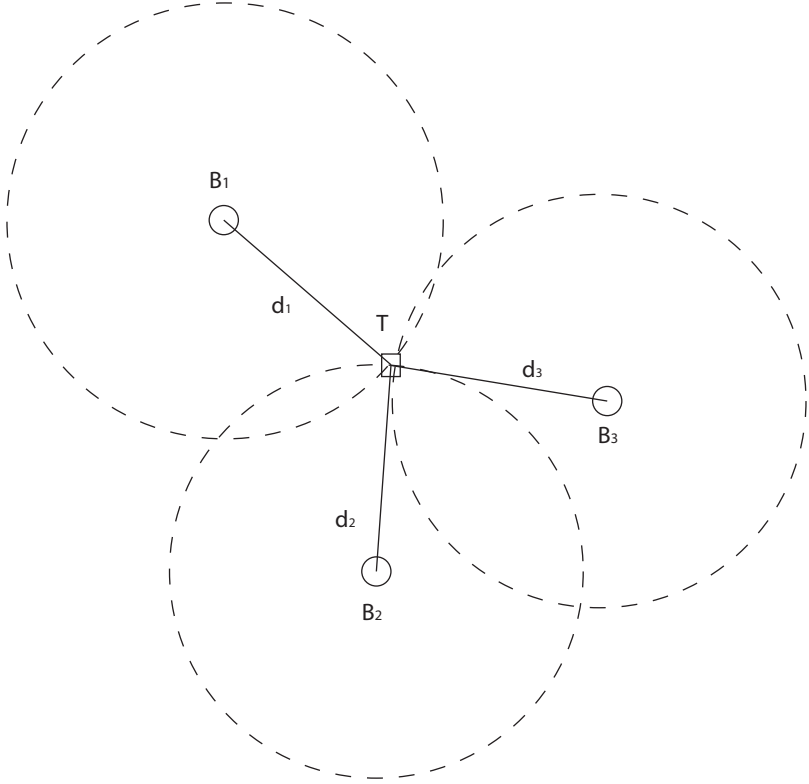}\caption{The dashed circles are the distances between the base stations $B_{i}$
and object $T$. The object $T$ is located at the intersection point
between the three dashed circles.\label{fig:intro} }
\end{figure}

Figure \ref{fig:intro} shows three base stations $B_{i}$ at known
positions $(a_{i},b_{i},c_{i})$, and one object $T$ at unknown position
$(x,y,z)$. The distances measurements $d_{i}$ between base stations
$B_{i}$ and object $T$ are known. The unkown position of object
$T$ can be estimated by the known positions of the base stations
$B_{i}$ and the distance measurements $d_{i}$. Measurement errors
are neglected in this paper, therefore distances measurements can
be referred as distances.

\subsection{Mathematical formulation\label{subsec:Mathematical-formulation}}

If the distance measurements between object $T$ and base stations
$B_{i}$ have no errors, then the unknown position of object $T$
 can be found by solving eq. (\ref{eq:obj1}) or eq. (\ref{eq:obj2}).
The distances between $B_{i}$ and $T$ are defined as

\begin{center}
\begin{tabular}{ccc}
$d_{i}^{2}=(x_{G}-a_{i})^{2}+(y_{G}-b_{i})^{2}+(z_{G}-c_{i})^{2}.$ &  & $1\leq i\leq N$\tabularnewline
\end{tabular}
\par\end{center}
\begin{itemize}
\item Objective function one:
\end{itemize}
\begin{equation}
F_{1}(x,y,z)\coloneqq\frac{1}{4}\sum_{i=1}^{N}\left[\sqrt{(x-a_{i})^{2}+(y-b_{i})^{2}+(z-c_{i})^{2}}-d_{i}\right]^{2}\label{eq:obj1}
\end{equation}
\begin{itemize}
\item Objective function two:
\end{itemize}
\begin{equation}
F_{2}(x,y,z)\coloneqq\frac{1}{4}\sum_{i=1}^{N}\left[(x-a_{i})^{2}+(y-b_{i})^{2}+(z-c_{i})^{2}-d_{i}^{2}\right]^{2}\label{eq:obj2}
\end{equation}

\cite{GPS1,GPS2}.

The solving of eq.(\ref{eq:obj1}) or eq.(\ref{eq:obj2}) can be done
by non-convex optimization \cite{LeastSquares} $F_{i}(x,y,z)\rightarrow argmin$.
Alternatively, the non-linear system can be transformed into a linear
system \cite{ION_JS,Sidorenko2016}. With the assumptions made in
the Section \ref{subsec:Mathematical-formulation} it is possible
to obtain a linear system. In more complex cases where the positions
of base stations $B_{i}$ are unknown this is not possible at all.
With regard to future extensions to determining the base station positions
as well as the location of the object $T$, this article focuses on
finding a solution with a non-convex optimization algorithm.

\subsection{Reason for the approach}

The objective functions (\ref{eq:obj1}) and (\ref{eq:obj2}) are
non-linear and non-convex. The optimization of the objective functions
can cause convergence to a local minimum $L$ instead of the global
minimum $G$ (see Table \ref{tab:Exmplation-of-the}). In our approach
instead the $F_{1}$ and $F_{2}$ the improved objective functions
$F_{L1}$ and $F_{L2}$ are used. Both have an additional variable
$\lambda$ compared to the $F$ functions.
\begin{itemize}
\item Improved objective function one:
\end{itemize}
\begin{equation}
F_{L1}(x,y,z,\lambda)\coloneqq\frac{1}{4}\sum_{i=1}^{N}\left[\sqrt{(x-a_{i})^{2}+(y-b_{i})^{2}+(z-c_{i})^{2}+\lambda^{2}}-d_{i}\right]^{2}
\end{equation}
\begin{itemize}
\item Improved objective function two:
\end{itemize}
\begin{equation}
F_{L2}(x,y,z,\lambda)\coloneqq\frac{1}{4}\sum_{i=1}^{N}\left[(x-a_{i})^{2}+(y-b_{i})^{2}+(z-c_{i})^{2}+\lambda^{2}-d_{i}^{2}\right]^{2}
\end{equation}

In the next section, we prove that the $F_{L2}$ has a saddle point
at every position of the local minimum $L(x_{L},y_{L},z_{L})$ of
$F_{2}$ . Therefore, the Levenberg-Marquardt algorithm has a lower
probability to converge to local minimum. 

\subsection{Characteristics of a local minimum \label{subsec:Assumptions-about-the}}

\subsubsection*{First premise }

The objective function has an unique global minimum at $G(x_{G},y_{G},z_{G})$
and at least one local minimum at $L(x_{L},y_{L},z_{L})$. \label{par:Assumption-1.:}

\subsubsection*{Second premise \label{subsec:Assumption-2. Main Prove}}

The first derivative of the $F_{L}$ with respect to x, y and z is
zero at the local minimum. Second derivative of the $F_{L}$ at the
same position is greater than zero (Table \ref{tab:Assumption-2}).

\begin{table}[H]
\caption{Premise two\label{tab:Assumption-2}}
\ \\
\begin{centering}
\begin{tabular}{c|c}
First derivative & Second derivative\tabularnewline
\hline 
$\left(\frac{\partial}{\partial x}F_{L}\right)(x_{L},y_{L},z_{L},0)=0$ & $\left(\frac{\partial^{2}}{\partial x^{2}}F_{L}\right)(x_{L},y_{L},z_{L},0)>0$\tabularnewline
 & \tabularnewline
$\left(\frac{\partial}{\partial y}F_{L}\right)(x_{L},y_{L},z_{L},0)=0$ & $\left(\frac{\partial^{2}}{\partial y^{2}}F_{L}\right)(x_{L},y_{L},z_{L},0)>0$\tabularnewline
 & \tabularnewline
$\left(\frac{\partial}{\partial z}F_{L}\right)(x_{L},y_{L},z_{L},0)=0$ & $\left(\frac{\partial^{2}}{\partial z^{2}}F_{L}\right)(x_{L},y_{L},z_{L},0)>0$\tabularnewline
\end{tabular}
\par\end{centering}
\end{table}

\subsection{Hypothesis}

The first derivative of the $F_{L}$, with respect to the additional
variable $\lambda$ is zero and the second derivative is less than
zero at the local minimum (See Table \ref{tab:Assumption-3}). In
combination with the first and second premise the local minimum becomes
a saddle point. The Levenberg-Marquardt (derivative based optimization
algorithm) would not converge to a saddle point.

\begin{table}[H]
\caption{Hypothesis \label{tab:Assumption-3}}
\ \\
\centering{}%
\begin{tabular}{c|c}
First derivative & Second derivative\tabularnewline
\hline 
$\left(\frac{\partial}{\partial\lambda}F_{L}\right)(x_{L},y_{L},z_{L},0)=0$ & $\left(\frac{\partial^{2}}{\partial\lambda^{2}}F_{L}\right)(x_{L},y_{L},z_{L},0)<0$\tabularnewline
 & \tabularnewline
 & $\left(\frac{\partial^{2}}{\partial x\partial\lambda}F_{L}\right)(x_{L},y_{L},z_{L},0)=0$\tabularnewline
 & \tabularnewline
 & $\left(\frac{\partial^{2}}{\partial y\partial\lambda}F_{L}\right)(x_{L},y_{L},z_{L},0)=0$\tabularnewline
 & \tabularnewline
 & $\left(\frac{\partial^{2}}{\partial z\partial\lambda}F_{L}\right)(x_{L},y_{L},z_{L},0)=0$\tabularnewline
\end{tabular}
\end{table}

\subsubsection*{Hypothesis for function\textmd{ $F_{1}$} }

Every local minima of function $F_{1}$ becomes a saddle point at
the same coordinates with function $F_{L1}$. We have no analytical
proof of this hypothesis, but the numerical results in Section \ref{sec:Numerical-results}
demonstrate its validity in practice.

\subsubsection*{Hypothesis for function\textmd{ $F_{2}$} }

Every local minima of function $F_{2}$ becomes a saddle point at
the same coordinates with function $F_{L2}$. This will be proven
in the following sections as well as demonstrated numerically.

\subsection{Proof of the hypothesis for objective function $F_{2}$\label{Appendix:Prove}}

In Section \ref{subsec:Assumptions-about-the} premises and the hypothesis
of our approach were introduced. In this section the hypothesis will
be proven for the objective function $F_{2}$. The proof of the hypothesis
for objective function $F_{1}$ has not been found yet. The empirical
results show that the approach works for both objective functions.
First, a new coordinate system is defined. This coordinate system
is centered in $L$ (local minimum) with $G$ (global minimum) on
the positive $x$ axis.

\subsubsection{Definition of the new coordinate system \label{subsec:CoordinateSystem}}

Without loss of generality the following coordinate system can be
used

\begin{center}
\begin{tabular}{|c|}
\hline 
$y_{G}=z_{G}=x_{L}=y_{L}=z_{L}=0$\tabularnewline
\hline 
\end{tabular}
\par\end{center}

and 

\begin{center}
\begin{tabular}{|c|}
\hline 
$x_{G}>0$.\tabularnewline
\hline 
\end{tabular}
\par\end{center}

The distances between the base stations $B_{i}$ and and object $T$
are 

\begin{equation}
d_{i}^{2}=(x_{G}-a_{i})^{2}+(b_{i})^{2}+(c_{i})^{2}.
\end{equation}

\subsubsection{First part of the auxiliary results }

The second objective function can also be written as

\begin{equation}
F_{2}(x,y,z)=\sum_{i=1}^{N}\varphi_{i}(x,y,z)^{2}
\end{equation}

with the auxiliary function $\varphi_{i}(x,y,z)$

\begin{equation}
\varphi_{i}(x,y,z):=(x-a_{i})^{2}+(y-b_{i})^{2}+(z-c_{i})^{2}-d_{i}^{2}.
\end{equation}

At the position of the local minimum $L$, the auxiliary function
becomes

\[
\varphi_{i}(0,0,0)=\left[a_{i}{}^{2}+b_{i}{}^{2}+c_{i}{}^{2}-d_{i}^{2}\right]=
\]

\[
=\left[a_{i}{}^{2}+b_{i}{}^{2}+c_{i}{}^{2}-(x_{G}-a_{i})^{2}-b_{i}{}^{2}-c_{i}{}^{2}\right]=
\]

\begin{equation}
=\left[a_{i}{}^{2}-(x_{G}-a_{i})^{2}\right]=\left[2a_{i}x_{G}-(x_{G})^{2}\right]=x_{G}\left[2a_{i}-x_{G}\right].\label{eq:aux1}
\end{equation}

Therefore, the second objective function at the local minimum can
be written as

\begin{equation}
F_{2}(0,0,0)=x_{G}{}^{2}\sum_{i=1}^{N}\left[2a_{i}-x_{G}\right]^{2}.\label{eq:OF_localMinima}
\end{equation}

\subsubsection{Simplification of the hypothesis}

In Section \ref{subsec:Assumptions-about-the} the premises for the
approach were presented. 
\begin{equation}
\left(\frac{\partial^{2}}{\partial\lambda^{2}}F_{L}\right)(x_{L},y_{L},z_{L},0)<0\label{eq:hypothesis}
\end{equation}

In the following it will be shown, that the hypothesis eq. (\ref{eq:hypothesis})
is always correct for the improved objective function two.

\begin{equation}
F_{L2}(x,y,z,\lambda)\coloneqq\frac{1}{4}\sum_{i=1}^{N}\left[(x-a_{i})^{2}+(y-b_{i})^{2}+(z-c_{i})^{2}+\lambda^{2}-d_{i}^{2}\right]^{2}
\end{equation}

Eq. (\ref{eqfirstDiffLamda}) and Eq. (\ref{eq:secondDiff}) are the
first and second derivative of objective function $F_{L2}$with respect
to $\lambda$.

\begin{equation}
\left(\frac{\partial}{\partial\lambda}F_{L2}\right)(x,y,z,\lambda)=\sum_{i=1}^{N}\left[(x-a_{i})^{2}+(y-b_{i})^{2}+(z-c_{i})^{2}+\lambda^{2}-d_{i}^{2}\right]\lambda\label{eqfirstDiffLamda}
\end{equation}

\begin{equation}
\left(\frac{\partial^{2}}{\partial\lambda^{2}}F_{L2}\right)(x,y,z,\lambda)=\sum_{i=1}^{N}\left[(x-a_{i})^{2}+(y-b_{i})^{2}+(z-c_{i})^{2}+\lambda^{2}-d_{i}^{2}\right]+2N\lambda^{2}\label{eq:secondDiff}
\end{equation}

At the local minimum $L(x_{L},y_{L},z_{L})$.

\[
\left(\frac{\partial^{2}}{\partial\lambda^{2}}F_{L2}\right)(0,0,0,0)=\sum_{i=1}^{N}\left[a_{i}{}^{2}+b_{i}{}^{2}+c_{i}{}^{2}-d_{i}^{2}\right]=\sum_{i}^{N}\varphi_{i}(0,0,0)=
\]
\begin{equation}
=x_{G}\sum_{i=1}^{N}\left[2a_{i}-x_{G}\right]=2x_{G}\sum_{i=1}^{N}a_{i}-N\,x_{G}{}^{2}.
\end{equation}

We want to show that $\left(\frac{\partial^{2}}{\partial^{2}\lambda}F_{L2}\right)(x_{L},y_{L},z_{L},0)<0$,
hence we have to prove the inequality eq. (\ref{eq:toProve}).

\begin{equation}
2x_{G}\sum_{i=1}^{N}a_{i}-N\,x_{G}{}^{2}<0
\end{equation}

\begin{equation}
2\sum_{i=1}^{N}a_{i}<N\,x_{G}\label{eq:toProve}
\end{equation}

\subsubsection{Second part of the auxiliary results }

In the next step the condition at the local minimum is analyzed. 

\begin{equation}
F_{2}(x,y,z)\coloneqq\frac{1}{4}\sum_{i=1}^{N}\left[(x-a_{i})^{2}+(y-b_{i})^{2}+(z-c_{i})^{2}-d_{i}^{2}\right]^{2}
\end{equation}

The first derivative of objective function two equates eq. (\ref{eq:fistDiff_x}), 

\[
\left(\frac{\partial}{\partial x}F_{2}\right)(x,y,z)=\sum_{i=1}^{N}\left[(x-a_{i})^{2}+(y-b_{i})^{2}+(z-c_{i})^{2}-d_{i}^{2}\right](x-a_{i})=
\]

\begin{equation}
=\sum_{i}^{N}\varphi_{i}(x,y,z)(x-a_{i})\label{eq:fistDiff_x}
\end{equation}

in combination with eq.(\ref{eq:aux1}) the first derivative becomes
eq. (\ref{eq: becomesTo}).

\[
\left(\frac{\partial}{\partial x}F_{2}\right)(0,0,0)=\sum_{i=1}^{N}\varphi_{i}(0,0,0)(-a_{i})=\sum_{i=1}^{N}x_{G}\left[2a_{i}-x_{G}\right](-a_{i})=
\]

\begin{equation}
=\left[x_{G}{}^{2}\sum_{i=1}^{N}a_{i}-2x_{G}\sum_{i=1}^{N}a_{i}{}^{2}\right].\label{eq: becomesTo}
\end{equation}

At the local minimum $L(x_{L},y_{L},z_{L})$ the first derivative
of objective function two equates zero. 

\begin{equation}
x_{G}{}^{2}\sum_{i=1}^{N}a_{i}-2x_{G}\sum_{i=1}^{N}a_{i}{}^{2}=0
\end{equation}

\begin{equation}
x_{G}\sum_{i=1}^{N}a_{i}=2\sum_{i=1}^{N}a_{i}{}^{2}\label{eq:aux2}
\end{equation}

This leads to that $\sum_{i=1}^{N}a_{i}>0.$ 

\subsubsection{Third part of the auxiliary results }

The objective function $F_{2}$ has a higher result at the local minimum
compared to the global minimum. It is assumed that the objective functions
have no errors, therefore the result of $F_{2}$ at the global minimum
has to be zero.

\begin{equation}
F_{2}(0,0,0)>F_{2}(x_{G},0,0)=0\label{eq:assumption}
\end{equation}

The term $F_{2}(0,0,0)$ of eq.(\ref{eq:assumption}) is replaced
by eq.(\ref{eq:OF_localMinima}). Eq. (\ref{eq:23}) can be converted
to eq. (\ref{eq:26}).

\begin{equation}
x_{G}{}^{2}\sum_{i=1}^{N}(2a_{i}-x_{G})^{2}>0\label{eq:23}
\end{equation}

\begin{equation}
\sum_{i=1}^{N}(2a_{i}-x_{G})^{2}>0
\end{equation}

\begin{equation}
4\sum_{i=1}^{N}a_{i}{}^{2}-4x_{G}\sum_{i=1}^{N}a_{i}+N\,x_{G}{}^{2}>0
\end{equation}

\begin{equation}
4x_{G}\sum_{i=1}^{N}a_{i}<4\sum_{i=1}^{N}a_{i}{}^{2}+N\,x_{G}{}^{2}\label{eq:26}
\end{equation}

In combination with eq.(\ref{eq:aux2}) the new inequality equates
eq. (\ref{eq:aux3}) .

\begin{equation}
8\sum_{i=1}^{N}a_{i}{}^{2}<4\sum_{i=1}^{N}a_{i}{}^{2}+N\,x_{G}{}^{2}
\end{equation}

\begin{equation}
4\sum_{i=1}^{N}a_{i}{}^{2}<N\,x_{G}{}^{2}\label{eq:aux3}
\end{equation}

\subsubsection{Proof by Cauchy-Bunyakovsky-Schwarz inequality}

The final step of the prove for the hypothesis, requires the Cauchy-Bunyakovsky-Schwarz
inequality \cite{SchwarzUngleichung} for $\mathbb{R}^{N}$.

\begin{center}
\begin{tabular}{cc}
What we want to prove that: & $2\sum_{i=1}^{N}a_{i}<N\,x_{G}$\tabularnewline
\end{tabular}
\par\end{center}

The Cauchy-Bunyakovsky-Schwarz inequality dictates that $\left|\left\langle \vec{x},\vec{y}\right\rangle \right|\leq\left\Vert \vec{x}\right\Vert \cdot\left\Vert \vec{y}\right\Vert $.
In our case the vectors are.

\begin{center}
\begin{tabular}{c}
$\vec{x}=\left(\begin{array}{c}
1\\
\vdots\\
1
\end{array}\right)$ and $\vec{y}=\left(\begin{array}{c}
a_{1}\\
\vdots\\
a_{n}
\end{array}\right)$\tabularnewline
\end{tabular}
\par\end{center}

The left term $2\sum_{i=1}^{N}a_{i}$ of eq.(\ref{eq:toProve}) is
due to the Cauchy-Bunyakovsky-Schwarz inequality smaller or equal
to $2\sqrt{N}\sqrt{\sum_{i=1}^{N}a_{i}{}^{2}}$.

\begin{equation}
2\sum_{i=1}^{N}a_{i}\leq2\sqrt{N}\sqrt{\sum_{i=1}^{N}a_{i}{}^{2}}\label{eq:Cauchy}
\end{equation}
\ \\

From eq.(\ref{eq:aux3}) it is known that $\sum_{i=1}^{N}(a_{i})^{2}<\frac{1}{4}N\cdot(x_{G})^{2}$,
therefore the right side of the eq.(\ref{eq:Cauchy}) can be written
as $2\sqrt{N}\sqrt{\frac{1}{4}N\,x_{G}{}^{2}}$.

The inequality becomes 
\begin{equation}
2\sum_{i=1}^{N}a_{i}<2\sqrt{N}\sqrt{\frac{1}{4}N\,x_{G}{}^{2}}=N\,x_{G}.
\end{equation}

\subsection{The effect of an additional variable on the global minimum}

The second derivative of $F_{L2}$with respect to $\lambda$ at the
global minimum is:

\begin{equation}
\left(\frac{\partial^{2}}{\partial\lambda^{2}}F_{L2}\right)(x_{G},y_{G},z_{G},\lambda_{G})=3\lambda_{G}^{2}N=0\,.
\end{equation}

At the global minimum the additional variable $\lambda$ has to be
zero and the second derivative must be bigger than zero. If the second
derivative is zero, a higher order derivative is required. 

\begin{equation}
\left(\frac{\partial^{3}}{\partial\lambda^{3}}F_{L2}\right)(x_{G},y_{G},z_{G},\lambda_{G})=\sum_{i=1}^{N}6\lambda_{G}N=0
\end{equation}

The third derivative is zero as well. Finally, the fourth derivative
is greater than zero, hence the additional variable has no effect
on the global minimum.

\begin{equation}
\left(\frac{\partial^{4}}{\partial\lambda^{4}}F_{L2}\right)(x_{G},y_{G},z_{G},\lambda_{G})=6N
\end{equation}

\subsection{No new local minima for $F_{L2}$ with $\lambda\protect\neq0$ \label{subsec:Local-minima-ofLOF}}

We have shown that the modified objective function $F_{L2}$ turns
the local minima of $F_{2}$ into saddle points and leaves the global
minimum unaffected. It remains to be proven that $F_{L2}$ does not
introduce new local minima that might adversely affect convergence
to the global minimum.

In this section we will show that in practically relevant base station
arrangements, $F_{L2}$ has no stationary points for $\lambda\neq0$
and $\mathbf{x}\neq\mathbf{x}_{G}$, and therefore no minima that
would lead an optimization method astray. We will show that if the
first derivative of $F_{L2}$ with respect to $\lambda$ vanishes
where $\lambda\neq0$, its gradient in the spacial directions is non-zero
for $\mathbf{x}\neq\mathbf{x}_{G}$. This proof is best presented
in vectorial notation. We will use $\mathbf{x}=(x,y,z)^{T}$ for the
position argument and $\mathbf{a}_{i}=(a_{i},b_{i},c_{i})^{T}$ for
the base station locations. 
\begin{eqnarray}
\frac{\partial}{\partial\lambda}\,F_{L2}(\mathbf{x},\lambda) & = & \;\lambda\sum\limits _{i}\left((\mathbf{x}-\mathbf{a}_{i})^{2}+\lambda^{2}-d_{i}^{2}\right)=0,\qquad\lambda\neq0\nonumber \\
 & \Rightarrow & \quad\sum\limits _{i}\left((\mathbf{x}-\mathbf{a}_{i})^{2}+\lambda^{2}-d_{i}^{2}\right)=0\label{eq:newminlambdasum}\\
\mbox{grad}_{\mathbf{x}}\,F_{L2}(\mathbf{x},\lambda) & = & \sum\limits _{i}\left((\mathbf{x}-\mathbf{a}_{i})^{2}+\lambda^{2}-d_{i}^{2}\right)\,(\mathbf{x}-\mathbf{a}_{i})\label{eq:grad}
\end{eqnarray}
Eq.~(\ref{eq:newminlambdasum}) allows us to add or subtract any
term not dependent on the summation index $i$ in the right-hand factor
of~(\ref{eq:grad}). We subtract $\mathbf{x}$ and add $\mathbf{a}_{\ast}=\frac{1}{N}\sum_{i=1}^{N}\mathbf{a}_{i}$,
the geometrical center of the base stations:

\begin{eqnarray*}
\mbox{grad}_{\mathbf{x}}\,F_{L2}(\mathbf{x},\lambda) & = & \sum\limits _{i}\left((\mathbf{x}-\mathbf{a}_{i})^{2}+\lambda^{2}-d_{i}^{2}\right)\,(\mathbf{x}-\mathbf{a}_{i}-\mathbf{x}+\mathbf{a}_{\ast})\\
 & = & -\sum\limits _{i}\left((\mathbf{x}-\mathbf{a}_{i})^{2}+\lambda^{2}-d_{i}^{2}\right)\,(\mathbf{a}_{i}-\mathbf{a}_{\ast})\,.
\end{eqnarray*}
By the construction of $\mathbf{a}_{\ast}$, we have $\sum_{i=1}^{N}(\mathbf{a}_{i}-\mathbf{a}_{\ast})=0$,
so now we can add or subtract any term not depending on the summation
index in the left-hand factor. We add $-\lambda^{2}-\mathbf{x}^{2}+\mathbf{x}_{G}^{2}$
and substitute $d_{i}=|\mathbf{x}_{G}-\mathbf{a}_{i}|$, expand the
squares and simplify, obtaining:

\begin{eqnarray*}
\mbox{grad}_{\mathbf{x}}\,F_{L2}(\mathbf{x},\lambda) & = & -\sum\limits _{i}\left((\mathbf{x}-\mathbf{a}_{i})^{2}-\mathbf{x}^{2}+\mathbf{x}_{G}^{2}-d_{i}^{2}\right)\,(\mathbf{a}_{i}-\mathbf{a}_{\ast})\\
 & = & -\sum\limits _{i}\left((\mathbf{x}-\mathbf{a}_{i})^{2}-\mathbf{x}^{2}+\mathbf{x}_{G}^{2}-(\mathbf{x}_{G}-\mathbf{a}_{i})^{2}\right)\,(\mathbf{a}_{i}-\mathbf{a}_{\ast})\\
 & = & \sum\limits _{i}(2\,\mathbf{x}\,\mathbf{a}_{i}-2\,\mathbf{x}_{G}\,\mathbf{a}_{i})\,(\mathbf{a}_{i}-\mathbf{a}_{\ast})\\
 & = & 2\;(\mathbf{x}-\mathbf{x}_{G})^{T}\sum\limits _{i}\mathbf{a}_{i}\otimes(\mathbf{a}_{i}-\mathbf{a}_{\ast})\\
 & = & 2\;(\mathbf{x}-\mathbf{x}_{G})^{T}\;\mathbf{M}\,.
\end{eqnarray*}
Here, $\mathbf{u}\otimes\mathbf{v}$ denotes the outer product, resulting
in a matrix with the entries $u_{i}\,v_{j}$. The matrix $\mathbf{M}$
can be brought into the following form: 
\begin{eqnarray*}
\mathbf{M} & = & \sum\limits _{i}\mathbf{a}_{i}\otimes\mathbf{a}_{i}-\left(\sum\limits _{i}\mathbf{a}_{i}\right)\otimes\mathbf{a}_{\ast}=\sum\limits _{i}\mathbf{a}_{i}\otimes\mathbf{a}_{i}-N\;\mathbf{a}_{\ast}\otimes\mathbf{a}_{\ast}\\
 & = & \sum\limits _{i}(\mathbf{a}_{i}-\mathbf{a}_{\ast})\otimes(\mathbf{a}_{i}-\mathbf{a}_{\ast})\,.
\end{eqnarray*}
The last step is analogous to the well-known derivation of the variance
of a data set. The result represents $\mathbf{M}$ as a sum of unnormalized
projection matrices onto the directions to the base stations from
their center.

By the calculation above, we have shown that the gradient of $F_{L2}$
has the form of a vector times a sum of projection matrices at all
local minima with $\lambda\neq0$. Projection matrices are positive
semidefinite by construction, and their null space is the sub-space
orthogonal to the projection direction. When adding several positive
semidefinite matrices, the null space of the result is the intersection
of the null spaces of the individual matrices, in our case the sub-space
orthogonal to all projection directions. For unambiguous location
in $n$ (2 or 3) dimensions, at least $n+1$ base stations are needed,
and they have to be arranged in a non-degenerate way, i.\ e.\ so
that the $\mathbf{a}_{i}-\mathbf{a}_{\ast}$ are a spanning set of
the whole space. This makes the matrix $\mathbf{M}$ positive definite,
and the gradient of $F_{L2}$ cannot be zero for $\mathbf{x}\neq\mathbf{x}_{G}$.
Therefore there are no local minima that prevent an optimization method
from converging to the global minimum.

\subsection{Two dimensional example }

In Section \ref{Appendix:Prove} it was proven that the $F_{L2}$
has a saddle point at the coordinates of the local minimum of $F_{2}$.
In this section an example is created with known coordinates of the
global $G(1,0)$ and local minimum $L(0,0)$. This example has the
aim to illustrate the converging steps of the Levenberg-Marquardt
algorithm for the $F_{2}$ and $F_{L2}$. The positions of the local
and global minimum leads to the coordinates of base stations $B_{i}$
(Table \ref{tab:Coordinates-of-objects} ). Further information can
be found in the appendix (\ref{sec:Possible-constellations-for}).

\begin{table}[H]
\caption{Coordinates of base stations $B_{i}$\label{tab:Coordinates-of-objects}}
\ \\
\centering{}%
\begin{tabular}{|c|c|c|}
\hline 
Base stations  & X-Position & Y-Position \tabularnewline
\hline 
\hline 
$B_{1}$ & 0 & 0\tabularnewline
\hline 
$B_{2}$ & $0.5$ & -2\tabularnewline
\hline 
$B_{3}$ & $0.5$ & 1\tabularnewline
\hline 
$B_{4}$ & $0.5$ & 3\tabularnewline
\hline 
\end{tabular}
\end{table}

\begin{figure}[H]
\centering{}\includegraphics[scale=0.6]{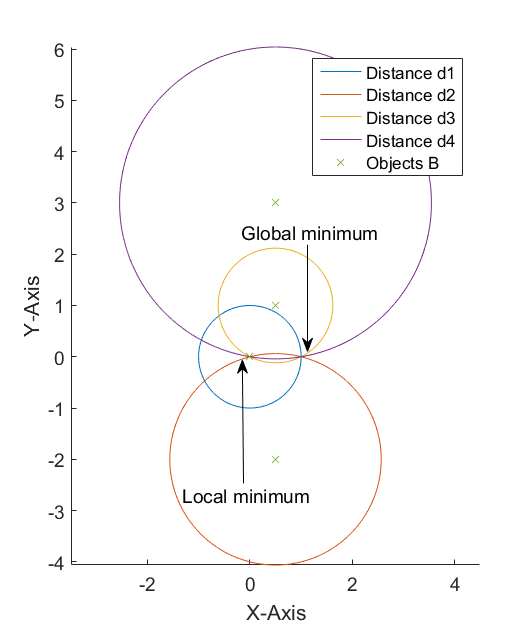}\caption{The circles represents the true distance between base stations $B_{i}$
and the global minimum. Blue circle is the distance between base station
$B_{1}$ and object $T$. Red circle is the distance between base
station $B_{2}$ and object $T$ . Yellow circle is the distance between
base station $B_{3}$ and object $T$ . Magenta circle is the distance
between base station $B_{4}$ and object $T$ respectively. \label{fig:Example-AB}}
\end{figure}

Figure \ref{fig:Example-AB}, shows the coordinates of the base stations
$B_{i}$, which are located in the center of the circles.

\begin{figure}[H]
\begin{centering}
\includegraphics[scale=0.4]{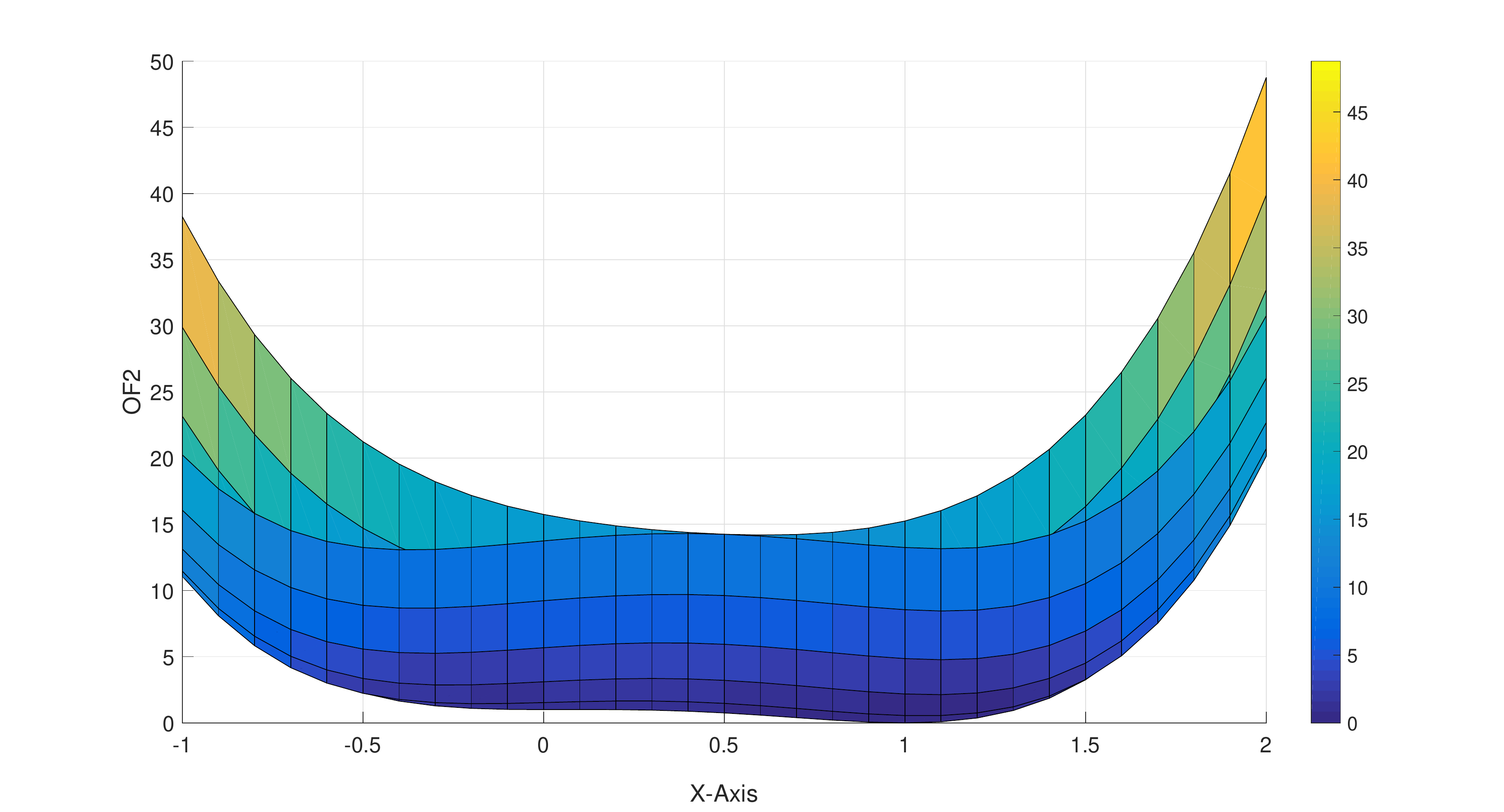}
\par\end{centering}
\begin{centering}
\includegraphics[scale=0.4]{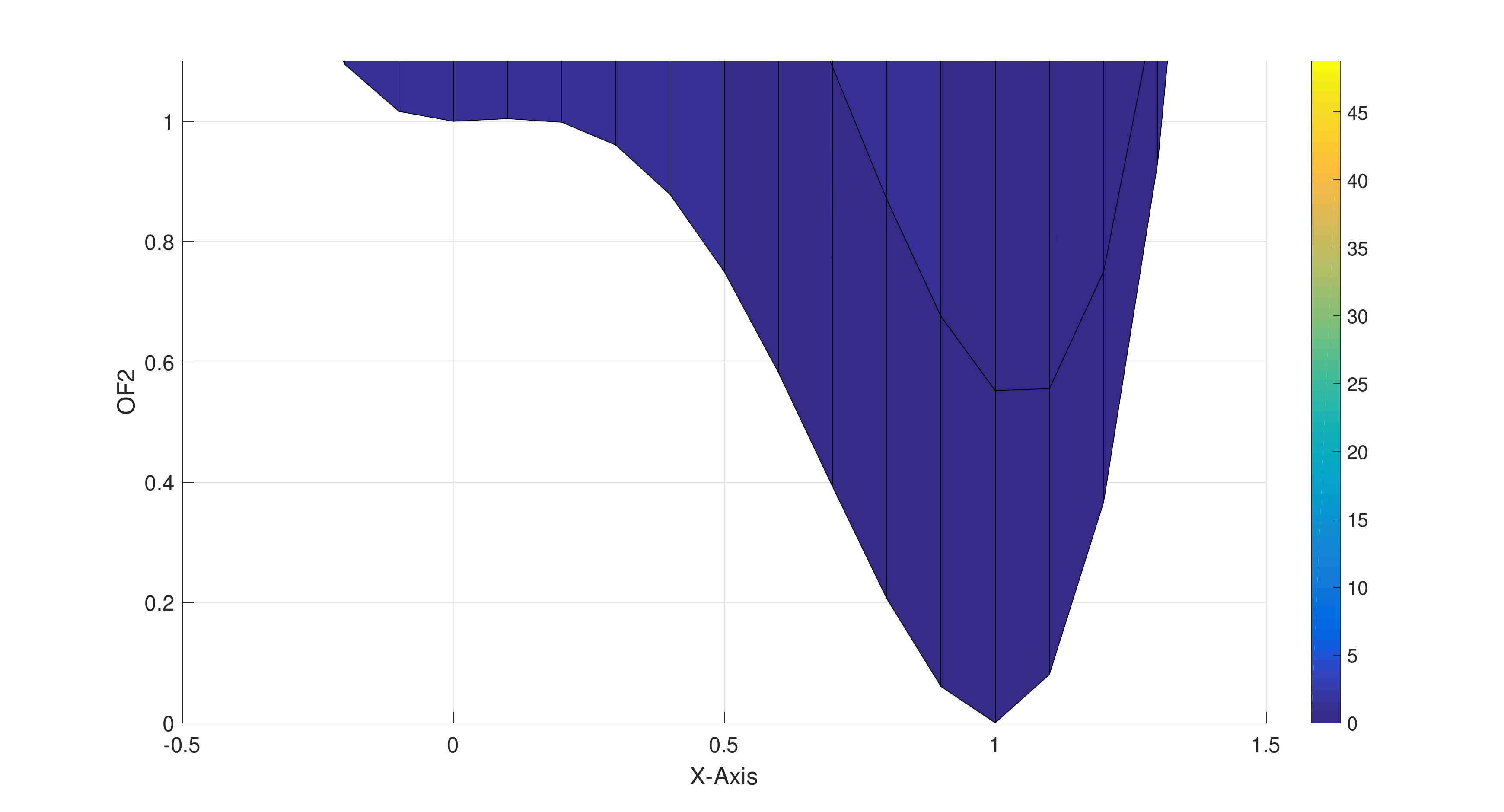}
\par\end{centering}
\centering{}\caption{Top Figure: Local minimum at $L(0,0)$ and global minimum at $G(1,0)$.
Colors from blue to yellow shows the residues of the objective function.
Bottom Figure is a zoom of top figure.\label{fig:Example-AB_sideView}}
\end{figure}

Figure \ref{fig:Example-AB_sideView}, shows the search space of objective
function $F_{2}$ and the zoom at the global minimum.

\subsubsection{Local optimization}

The Levenberg-Marquardt algorithm uses the derivative to obtain the
stepsize, therefore it is important that the initial estimate for
the additional variable $\lambda$ is non-zero. Otherwise $\lambda$
remains zero, and $F_{L2}$ is effectively reduced to $F_{2}$.

Table \ref{tab:Initial-position} shows initial estimates of the optimization.

\begin{table}[H]
\caption{Iteration steps of the Levenberg-Marquardt for $F_{2}$and $F_{L2}$.
\label{tab:Initial-position}}
\ \\
\centering{}%
\begin{tabular}{|c|c|c|c|}
\hline 
 & $x$ & $y$ & $\lambda$\tabularnewline
\hline 
\hline 
Initial estimate & -1 & 2 & 1\tabularnewline
\hline 
\end{tabular}
\end{table}

\begin{figure}[H]
\centering{}\includegraphics[scale=0.4]{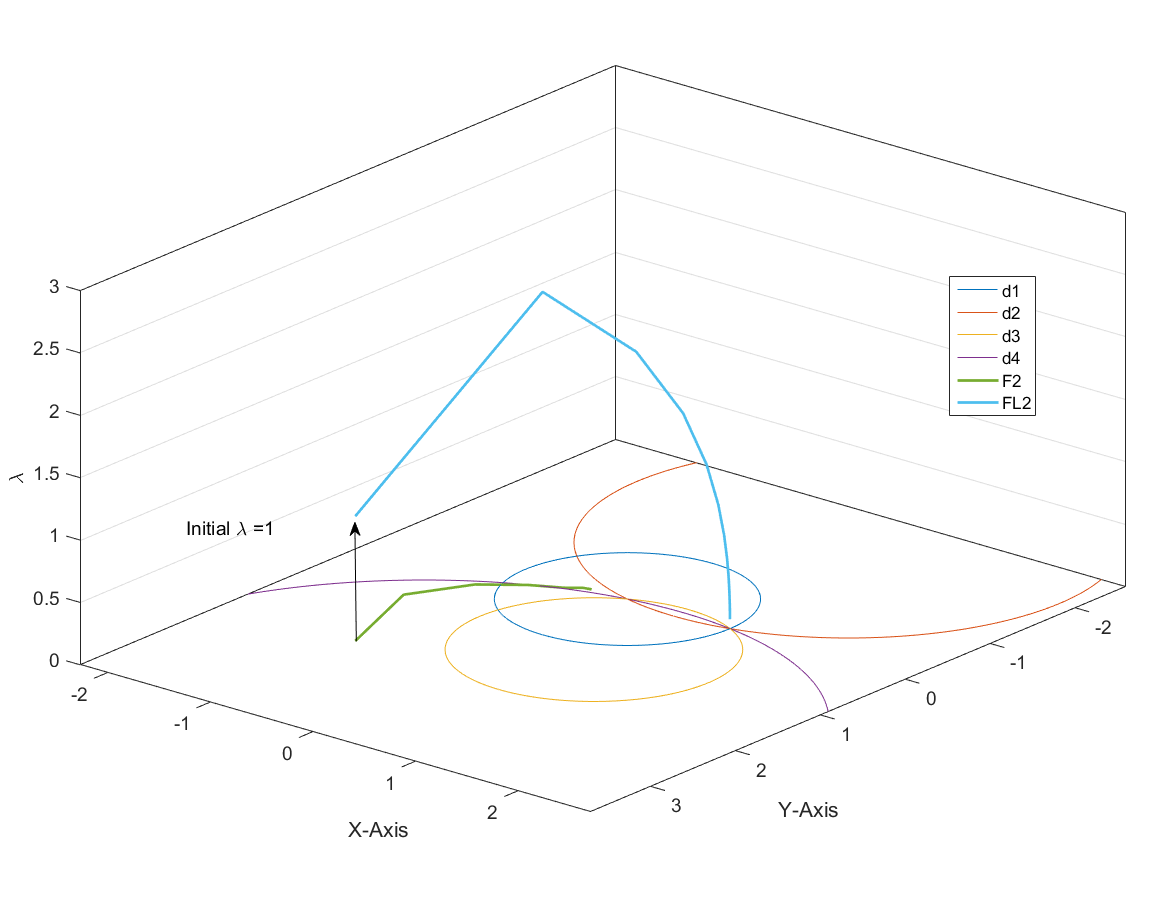}\caption{Iteration steps of the Levenberg-Marquardt for $F_{2}$and $F_{L2}$.
F2: Objective function $F_{2}$. FL2: Improved objective function
two $F_{L2}$. Blue line: Opimization steps of $F_{2}$. Gree line:
Optimization steps of $F_{L2}$. The circles blue, bed, yellow and
magenta are the distances between base stations $B_{i}$ and object
$T$ \label{fig:Optimization:-AB}}
\end{figure}

\begin{figure}[H]
\centering{}\includegraphics[scale=0.37]{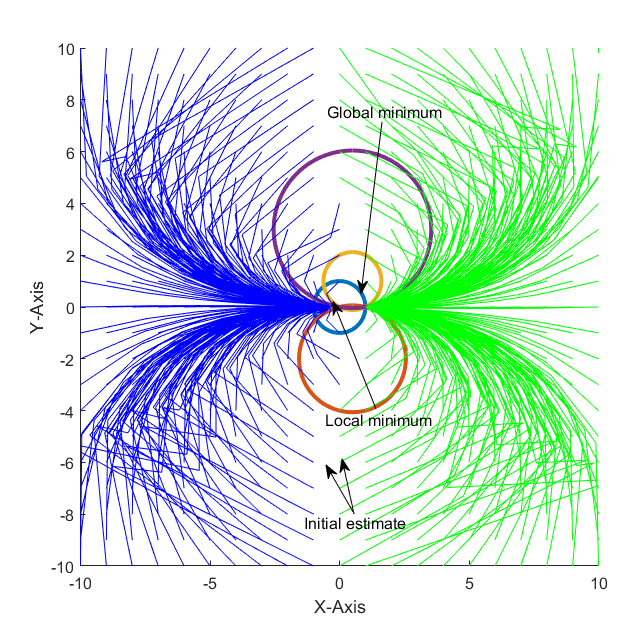}\includegraphics[scale=0.37]{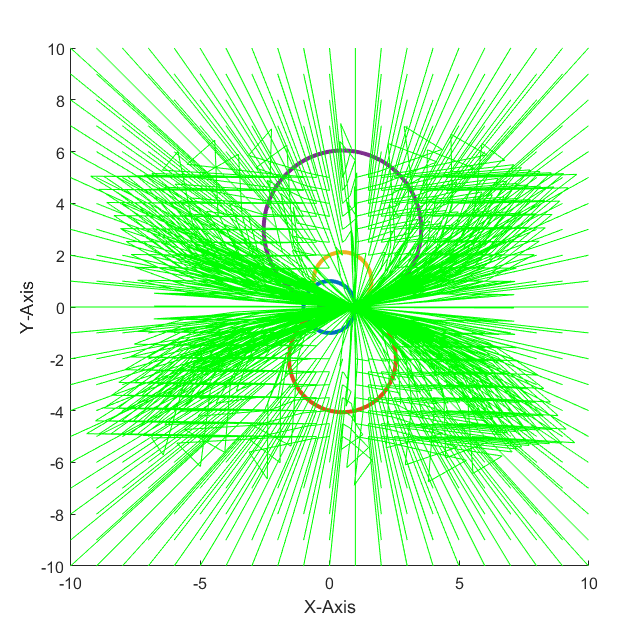}\caption{Left figure shows $F_{2}$ and the right figure shows the $F_{L2}$
with different initial estimates. Green: Convergence to global minimum.
Blue: Convergence to local minimum. \label{fig:2D-konvergenz-radius}}
\end{figure}

In Figure \ref{fig:Optimization:-AB} the result of the optimization
can be observed. The blue path shows the steps of the improved objective
function $F_{L2}$, which converge to the global minimum $G(1,0)$.
On the other hand, the original objective function $F_{2}$ represented
by the green line, converges to the local minimum  $L(0,0)$. In Figure
\ref{fig:2D-konvergenz-radius}, the optimization was repeated with
different initial estimates. The blue lines begin at the initial estimate
and end at the local minimum $L(0,0)$. The green lines end at the
global minimum $G(1,0)$. In this constellation, every start estimate
with $x<0$ converges to the local minimum. The $F_{L2}$ always converges
to the global minimum $G(1,0)$.

\section{Numerical results\label{sec:Numerical-results}}

\begin{figure}[H]
\begin{centering}
\includegraphics[scale=0.5]{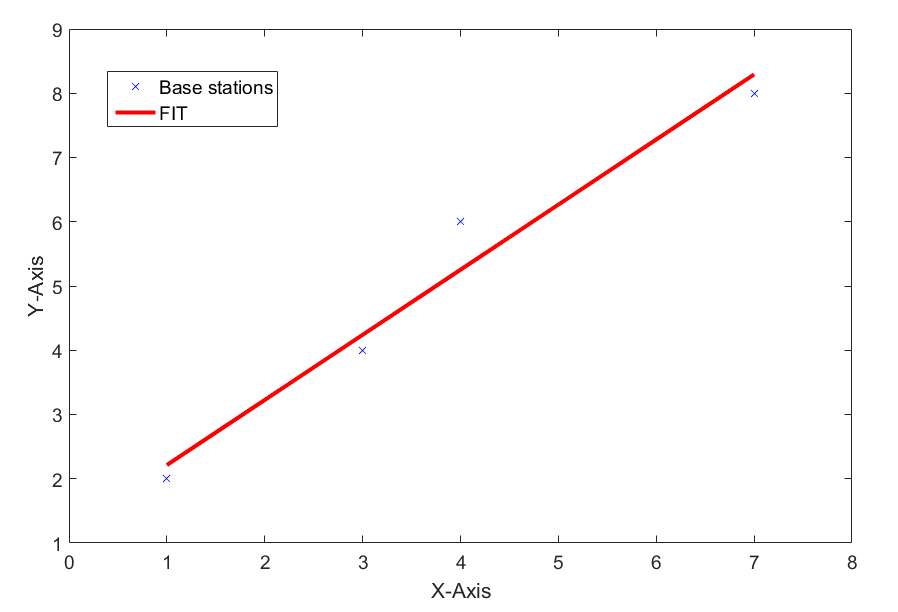}
\par\end{centering}
\caption{Close to collinearity. Red line is fit between the base stations $B$
. \label{fig:HDOP}}
\end{figure}

The tests were carried out with MALTAB Levenberg-Marquardt algorithm
at default settings (Table. \ref{tab:Default-Matlab-'Levenberg}).
\begin{table}[H]
\caption{Default MATLAB 'Levenberg Marquardt algorithm' parameter \label{tab:Default-Matlab-'Levenberg}}
\ \\
\centering{}{\footnotesize{}}%
\begin{tabular}{|c|c|}
\hline 
 & {\footnotesize{}Value}\tabularnewline
\hline 
\hline 
{\footnotesize{}Maximum change in variables for finite-difference
gradients} & {\footnotesize{}Inf}\tabularnewline
\hline 
{\footnotesize{}Minimum change in variables for finite-difference
gradients } & {\footnotesize{}0}\tabularnewline
\hline 
{\footnotesize{}Termination tolerance on the function value} & {\footnotesize{}1e-6}\tabularnewline
\hline 
{\footnotesize{}Maximum number of function evaluations allowed} & {\footnotesize{}100{*}numberOfVariables}\tabularnewline
\hline 
{\footnotesize{}Maximum number of iterations allowed} & {\footnotesize{}400}\tabularnewline
\hline 
{\footnotesize{}Termination tolerance on the first-order optimality} & {\footnotesize{}1e-4 }\tabularnewline
\hline 
{\footnotesize{}Termination tolerance on x} & {\footnotesize{}1e-6}\tabularnewline
\hline 
{\footnotesize{}Initial value of the Levenberg-Marquardt parameter} & {\footnotesize{}1e-2}\tabularnewline
\hline 
\end{tabular}{\footnotesize \par}
\end{table}
 The base stations $B_{i}$, object $T$ and initial estimates were
randomly generated in a 10x10x10 cube. Figure \ref{fig:HDOP} shows
an unfavorable constellation of the base stations close to collinearity. 

This constellations have been avoided by the requirement that every
normalized singular value of the covariance matrix should be higher
than $0.1$. 

\begin{itemize}
\item Error term:
\end{itemize}
\begin{equation}
E=\sum_{j=1}^{M}\sqrt{(x-x_{G})^{2}+(y-y_{G})^{2}+(z-z_{G})^{2}}\label{eq:errorTerm}
\end{equation}

\subsection{Results with the objective function $F_{1}$}

\begin{figure}[H]
\begin{centering}
\includegraphics[scale=0.5]{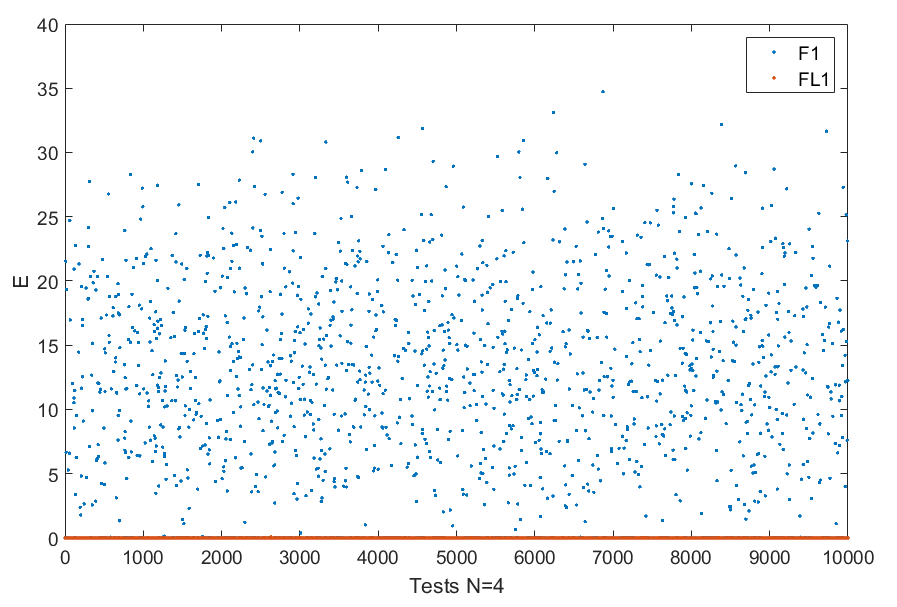}
\par\end{centering}
\caption{Blue dots: Objective function $F_{1}$. Red dots: Improved objective
function $F_{L1}$\label{fig:results1}}
\end{figure}

In the following section the results of the optimization with a two
dimensional $F_{1}$ are presented. Figure \ref{fig:results1} shows
the error term with different constellations of the four base stations
$B_{i}$. It can be seen that the $F_{L1}$ has no outlier. It has
not been proven yet, that the local minimum of $F_{1}$ becomes a
saddle point with $F_{L1}$ . However, the results show a significant
effect of the $F_{L1}$ on the optimization process.

\subsection{Results with the objective function $F_{2}$}

In the following section the results of the optimization with $F_{2}$
and $F_{L2}$ with four base stations $B_{i}$ are presented. In this
case it was proven that the local minimum of $F_{2}$ becomes a saddle
point with an additional variable. The error term eq. (\ref{eq:errorTerm})
of $F_{L2}$ is always smaller than 0.5 hence, it can be assumed that
the $F_{L2}$ has never converges to local minimum.

\begin{figure}[H]
\begin{centering}
\includegraphics[scale=0.5]{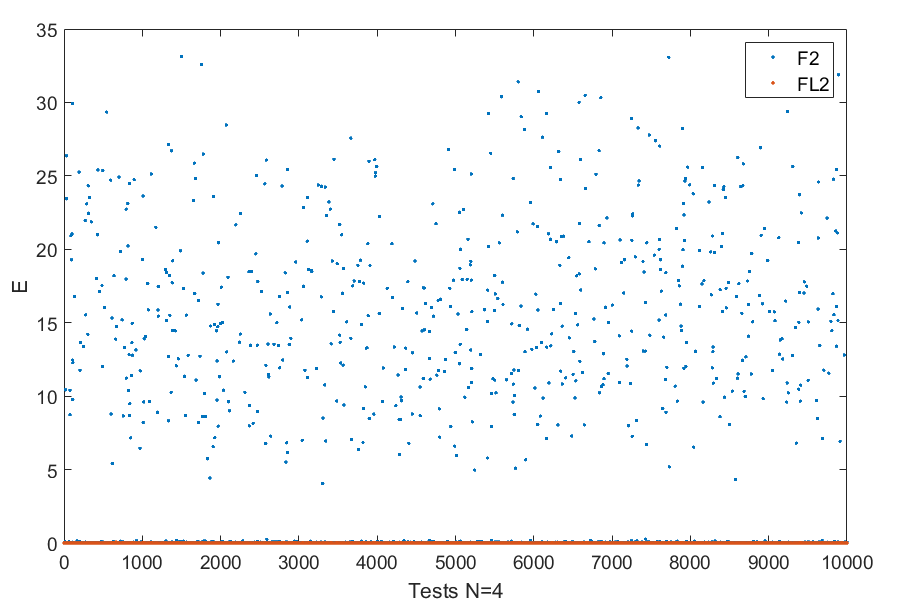}
\par\end{centering}
\caption{The blue dots are the results of the error term of $F_{2}$ . The
red dots are the results of the $F_{L2}$ }
\end{figure}

\subsection{Summary of the results}

Table \ref{tab:results} and \ref{tab:results3D}, shows the summary
of the obtained results. For every number of objects (N), 10.000 constellations
have been created and tested with Levenberg-Marquardt. The $F_{L}$
has not a single time convergences into a local minimum. 

\begin{table}[H]
\centering{}\caption{The examples are based on a 2-D model. N: Number of base stations
$B_{i}$. $F_{1}$: Objective function one, $F_{2}$: Objective function
two, M: Mean error, Sigma: Standard deviation, L: Amount of local
minima (Error bigger then $0.5$) \label{tab:results}. }
\ \\
\begin{tabular}{|c|c|c|c|c|c|}
\hline 
N & Objective function & $F_{1}$: $M\pm\sigma$ & $F_{1}$: L & $F_{2}$: $M\pm\sigma$ & $F_{2}$: L\tabularnewline
\hline 
\hline 
4 & $F$ & $1.8757\pm5.3120$ & 1357 & $1.0376\pm4.2345$ & 634\tabularnewline
\hline 
4 & $F_{L}$ & \textbf{$0.0015\pm0.0017$} & \textbf{0} & \textbf{$0.0020\pm0.0022$} & \textbf{0}\tabularnewline
\hline 
\hline 
5 & $F$ & $1.3981\pm4.6856$ & 982 & $0.6917\pm3.5277$ & 399\tabularnewline
\hline 
5 & $F_{L}$ & \textbf{$0.0014\pm0.0014$} & \textbf{0} & \textbf{$0.0019\pm0.0019$} & \textbf{0}\tabularnewline
\hline 
\hline 
6 & $F$ & $1.2150\pm4.4387$ & 810 & $0.5120\pm3.0919$ & 286\tabularnewline
\hline 
6 & $F_{L}$ & \textbf{$0.0014\pm0.0013$} & \textbf{0} & \textbf{$0.0019\pm0.0017$} & \textbf{0}\tabularnewline
\hline 
\hline 
7 & $F$ & $0.9128\pm3.9111$ & 586 & $0.3501\pm2.5899$ & 182\tabularnewline
\hline 
7 & $F_{L}$ & \textbf{$0.0014\pm0.0013$} & \textbf{0} & \textbf{$0.0018\pm0.0016$} & \textbf{0}\tabularnewline
\hline 
\end{tabular}\\
\end{table}

\begin{table}[H]
\centering{}\caption{The examples are based on a 3-D model. N: Number of base stations
$B_{i}$. $F_{1}$: Objective function one, $F_{2}$: Objective function
two, M: Mean error, Sigma: Standard deviation, L: Amount of local
minima (Error bigger then $0.5$) \label{tab:results3D}}
\ \\
\begin{tabular}{|c|c|c|c|c|c|}
\hline 
N & Objective function & $F_{1}$: $M\pm\sigma$ & $F_{1}$: L & $F_{2}$: $M\pm\sigma$ & $F_{2}$: L\tabularnewline
\hline 
\hline 
7 & $F$ & $0.6198\pm3.0085$ & 494 & $0.3230\pm2.2790$ & 216\tabularnewline
\hline 
7 & $F_{L}$ & \textbf{$0.0012\pm0.0011$} & \textbf{0} & \textbf{$0.0017\pm0.0014$} & \textbf{0}\tabularnewline
\hline 
\end{tabular}
\end{table}

\section{Discussion}

In the methodology section it was proven that the improved objective
function $F_{L2}$ has no local minima for non-trivial constellations.
This was underpinned by more than 100,000 numerical tests with reasonable
constellations. It has not been proven yet, that the local minimum
of $F_{1}$ becomes a saddle point with $F_{L1}$ . However, the results
show a significant effect of the $F_{L1}$ on the optimization process.
The objective function $F_{2}$ performed better than the objective
function $F_{1}$. Furthermore, was the number of false results $L$
decreasing with a higher amount of base stations $B_{i}$. It is important
that the initial estimate of the additional variable is unequal zero.
Otherwise gradient-based optimization algorithms like Levenberg-Marquardt
would not converge to the additional dimension. In all test scenarios
the positions of base stations $B_{i}$were known. Under the following
conditions it is also possible to obtain the solution analytically.
In the case of unknown positions of base stations $B_{i}$ and objects
$T_{j}$ it is not feasible anymore. At this point, our approach becomes
extremely valuable.

\section{Appendix}

\section{Possible constellations for the example \label{sec:Possible-constellations-for}}

The coordinate system was described in Section \ref{subsec:CoordinateSystem}.
We want to find base station constellations with a local minimum at
$x_{L}$,$y_{L}$,$z_{L}$ and a global minimum at $x_{G}$,$y_{G}$,$z_{G}$
.

\subsection{First derivative}

The first derivative of objective function two has to be zero at the
local minimum.

\begin{equation}
\left(\frac{\partial}{\partial x}F_{2}\right)(0,0)=\sum_{i=1}^{N}\left[(a_{i}{}^{2}-(x_{G}-a_{i})^{2})\cdot(-a_{i})\right]=0\label{eq:derivativeOneX}
\end{equation}

\begin{equation}
\left(\frac{\partial}{\partial y}F_{2}\right)(0,0)=\sum_{i=1}^{N}\left[(a_{i}{}^{2}-(x_{G}-a_{i})^{2})\cdot(-b_{i})\right]=0\label{eq:derivativeOneY}
\end{equation}

In Table \ref{tab:Obvious-solutions} obvious solutions of eq. (\ref{eq:derivativeOneX})
and eq. (\ref{eq:derivativeOneY}) can be found.

\begin{table}[H]
\caption{Obvious solutions \label{tab:Obvious-solutions}}
\ \\
\centering{}%
\begin{tabular}{c|cc}
Solution & X-Axis & Y-Axis\tabularnewline
\hline 
 &  & \tabularnewline
Solution A & $\forall i\,:\,a_{i}=0$ & $\forall i\,:\,b_{i}=0$\tabularnewline
 &  & \tabularnewline
Solution B & $\forall i\,:\,a_{i}=\frac{1}{2}x_{G}$ & $\mathbb{R}$\tabularnewline
\end{tabular}
\end{table}

With different combinations of solution A and B, it is possible to
create further solutions. 

\subsection{Second derivative}

The second derivative of objective function two with respect to x,
has to be higher than zero at the local minimum.

\begin{equation}
\left(\frac{\partial^{2}}{\partial x^{2}}F_{2}\right)(0,0)=\sum_{i=1}^{N}\left[2a_{i}^{2}+a_{i}^{2}-(x_{G}-a_{i})^{2}\right]>0\label{eq:derivative2}
\end{equation}

Insert eq.(\ref{eq:derivativeOneX}) into eq.(\ref{eq:derivative2})
leads to the first condition eq. (\ref{eq:condition1}).

\begin{equation}
\sum_{i=1}^{N}3a_{i}>Nx_{G}.\label{eq:condition1}
\end{equation}

The number of base stations $B_{i}$, which use solution B should
be higher than the number of objects with solution A. Otherwise there
will be no local minimum at the coordinates $L(0,0)$.

\begin{table}[H]
\caption{Notation of variable S \label{tab:Definition-of-variableS}}
\ \\
\centering{}%
\begin{tabular}{|c|c|}
\hline 
Notation & Number of base stations provided by\tabularnewline
\hline 
\hline 
$S_{1}$ & Solution A\tabularnewline
\hline 
$S_{2}$ & Solution B\tabularnewline
\hline 
\end{tabular}
\end{table}

Solution B cause to a positive second derivative $2\sum_{i=1}^{S_{2}}a_{i}{}^{2}$
and solution A to a negative $-x_{G}^{2}$. The second derivative
has to be higher than zero $\frac{\partial^{2}}{\partial x^{2}}F_{2}(0,0)>0$,
therefore $2\sum_{i=1}^{S_{2}}a_{i}{}^{2}>\sum_{i=1}^{S_{1}}x_{G}^{2}$.
The left term can be replaced by $\sum_{i=1}^{S_{2}}\frac{1}{4}x_{G}^{2}$.
This leads to,

\[
\frac{1}{2}S_{2}>S_{1}.
\]

The second derivative of objective function two with respect to y,
has to be positive at the local minimum $L(0,0)$.

\begin{equation}
\left(\frac{\partial^{2}}{\partial y^{2}}F_{2}\right)(0,0)=\sum_{i=1}^{N}\left[2b_{i}^{2}-x_{G}^{2}+2a_{i}x_{G}\right]>0\label{eq:derivative2Y}
\end{equation}

With solution B eq. (\ref{eq:derivative2Y}) becomes 

\[
\left(\frac{\partial^{2}}{\partial y^{2}}F_{2}\right)(0,0)=2\sum_{i=1}^{S_{2}}b_{i}^{2}.
\]

Solution A cause to a negative result of the second derivative, therefore

\[
2\sum_{i=1}^{S_{2}}b_{i}^{2}>x_{G}^{2}S_{1}.
\]

All the conditions for a local minimum at $L(0,0)$ can be found in
Table \ref{tab:Conditions-for-localMinima}.

\begin{table}[H]
\caption{Conditions for a local minimum at $L(0,0)$ \label{tab:Conditions-for-localMinima}}
\ \\
\centering{}%
\begin{tabular}{c|c}
 & Conditions\tabularnewline
\hline 
 & \tabularnewline
1 & $3\sum_{i=1}^{N}a_{i}>Nx_{G}$\tabularnewline
 & \tabularnewline
2 & $0.5\cdot S_{2}>S_{1}$\tabularnewline
 & \tabularnewline
3 & $2\sum_{i=1}^{S_{2}}b_{i}^{2}>x_{G}^{2}S_{1}$\tabularnewline
\end{tabular}
\end{table}

\subsection{Used constellations in the example}

The used assumptions for the example and the coordinates of the base
stations $B_{i}$ can be found in Table \ref{tab:Assumptions-used-for-Example}
and Table \ref{tab:Coordinates-of-objectB}.

\begin{table}[H]
\caption{Assumptions used for the example \label{tab:Assumptions-used-for-Example}}
\ \\
\centering{}%
\begin{tabular}{|c|}
\hline 
Assumptions\tabularnewline
\hline 
\hline 
$S_{1}=1$\tabularnewline
\hline 
$S_{2}=3$\tabularnewline
\hline 
\end{tabular}
\end{table}

\begin{table}[H]
\caption{Coordinates of object B \label{tab:Coordinates-of-objectB}}
\ \\
\centering{}%
\begin{tabular}{|c|c|c|}
\hline 
Base stations $B_{i}$with index & X-Axis  & Y-Axis\tabularnewline
\hline 
\hline 
1 & 0 & 0\tabularnewline
\hline 
2 & $0.5\cdot x_{G}$ & -2\tabularnewline
\hline 
3 & $0.5\cdot x_{G}$ & 1\tabularnewline
\hline 
4 & $0.5\cdot x_{G}$ & 3\tabularnewline
\hline 
\end{tabular}
\end{table}

\section*{Acknowledgments}

The first author would like to thank Sebastian Bullinger, Gregor Stachowiak,
Sebastian Tome for their inspiring discussion and the Fraunhofer IOSB
for making this work possible.

\bibliographystyle{plain}
\bibliography{reference}

\end{document}